\documentclass[12pt]{article}

\usepackage{amscd,amsmath,amssymb}
\newtheorem{definition}{Definition}
\newtheorem{theorem}{Theorem}
\newtheorem{lemma}{Lemma}
\newtheorem{proposition}{Proposition}
\newtheorem{corollary}{Corollary}

\makeatletter
\def\newproof#1{\@nprf{#1}}

\def\@nprf#1#2{\@xnprf{#1}{#2}}

\def\@xnprf#1#2{\expandafter\@ifdefinable\csname #1\endcsname
\global\@namedef{#1}{\@prf{#1}{#2}}\global\@namedef{end#1}{\@endproof}}

\def\@prf#1#2{\@xprf{#1}{#2}}

\def\@xprf#1#2{\@beginproof{#2}{\csname the#1\endcsname}\ignorespaces}

\def\@beginproof#1{\rm \trivlist \item[\hskip \labelsep{\it #1.\/}]}

\def\qed{\sqcap\!\!\!\!\sqcup}
\def\@endproof{\hfill$\qed$\endtrivlist}

\newproof{@proof}{Proof}
\newenvironment{proof}{\begin{@proof}}{\end{@proof}}

\makeatother

\def\de{\delta}
\def\De{\Delta}
\def\th{\theta}
\def\eps{\varepsilon}
\def\si{\sigma}
\def\Si{\Sigma}
\def\al{\alpha}
\def\Ga{\Gamma}
\def\ga{\gamma}
\def\La{\Lambda}
\def\la{\lambda}
\def\ti{\tilde}
\def\kappa{\varkappa}

\def\Bbb{\mathbb}

\title{Distinguished properties of the gamma process, and related topics}
\author{N.~Tsilevich%
\thanks{St.~Petersburg State University, St.~Petersburg, Russia.
Partially supported by RFBR grant 99--01--00098 and
DFG--RFBR grant 99--01--04027.
Email: {\tt natalia@pdmi.ras.ru}}
\and A.~Vershik
\thanks{Steklov Institute of Mathematics at St.~Petersburg, St.~Petersburg, Russia.
Partially supported by INTAS grant 93570-ext and
RFBR grant 99--01--00098.
Email: {\tt vershik@pdmi.ras.ru}} 
\and M.~Yor%
\thanks{Laboratoire de Probabilit\'es et Modeles al\'eatoires, Tour 56, 4 Place Jussieu,
75252 Paris Cedex 05, France}
}

\begin{document}
\maketitle

\begin{abstract}
We study fundamental properties of the gamma process and their relation to
       various topics such as Poisson--Dirichlet measures and stable processes.  
       We prove the quasi-invariance of the gamma process
       with respect to a large group of linear transformations. We also show
       that it is a renormalized limit of the stable
       processes and has an equivalent sigma-finite measure (quasi-Lebesgue) 
       with important invariance properties. New properties of the gamma process
       can be applied to the Poisson---Dirichlet measures. We also emphasize the deep 
       similarity  between the gamma process and the Brownian motion. 
       The connection of the above topics makes more transparent 
       some old and new facts about stable and gamma processes, and the 
       Poisson-Dirichlet measures.

\smallskip
{\it Keywords.} L\'evy processes, gamma process, stable processes, 
Po\-is\-son--Dirichlet measures, multiplicative quasi-invariance, 
quasi-Le\-bes\-gue measure, Markov--Krein identity. 
\end{abstract}

\tableofcontents

\section{Introduction}\label{s:intro}

The purpose of this work is to link various questions from the
probability theory, combinatorics, and representation theory, which
relate to the gamma process and the Poisson--Dirichlet measures, and which were
not regarded as a single whole. Since in spite of many aspects of our 
work its main object is the gamma process, first
we would like to present some fundamental properties of
the classical gamma process and its role for further considerations.

The standard gamma distribution (e.g.~see~\cite{JK70})
with parameter $\al>0$
is a distribution on the positive half-line with the following density:
$$ 
\frac{u^{\alpha-1}e^{-u}}{\Gamma(\alpha)},\quad u >0.
$$
Random variable with the gamma distribution is called a gamma variable.
The classical gamma process ($\gamma_t$, $t \geq 0$) is an increasing 
stochastic process with independent homogeneous increments such 
that $\ga_t$ has the gamma distribution with parameter $t$ for each $t>0$.

The following fundamental properties of this process
have led us to more general developments
and definitions discussed below, and we suggest that the reader 
comes back to them as a ``toy model'' when reading more general 
later parts of our paper.
            
\begin{description}
\item[1. Decomposition property:]
        If two variables $\xi$, $\phi$ are independent and have standard gamma 
        distributions, then $\xi + \phi$ and $\frac{\xi}{\xi+\phi}$ are 
        also independent.
        Moreover, this is a {\it characteristic property} of the gamma distribution
        (Lukacs' theorem, see Section~\ref{s:gamma} and~\cite{JK70}).
        Namely, if $x_1$, $x_2$
         are independent positive random variables, and $x_1+x_2$ and
         $\frac{x_1}{x_1+x_2}$ are also independent, then $x_1$ and $x_2$
         are common multiples of standard gamma variables with some
         parameters.
\end{description}
The first property easily implies 
\begin{description}
\item[2. Independence property:]         
         For each $t>0$, the process ($\gamma_u/{\gamma_t}$, $u \leq t$)
         is independent on ($\gamma_v$, $v \geq t$).
\end{description}
From property~2 one can deduce the following 
\begin{description}
\item[3. Quasi-invariance property:]
        For each $t>0$ and $a>0$ the law 
        of the process $(\frac{\gamma_u}{1+a},\, u \leq t)$ 
        is equivalent to that of              
         $(\gamma_u,\, u \leq t)$ with Radon--Nikodym density
$(1+a)^t \cdot \exp(-a\gamma_t)$.           
\end{description}
         This property will lead us to a very important quasi-invariance 
         of the law of the gamma process which was discovered in~\cite{GGV83} 
         in a different way and was used  for the representation theory 
         (see Section~\ref{s:main}).
\begin{description}
\item[4. Gamma distribution as a limit of stable distributions:]
       It was noticed that for a small index $\alpha$ the     
       stable distribution could be renormalized in such a way that the limit
       when $\alpha$ tends to zero is the gamma distribution~\cite{Cr75, ShS77}.
           This fact leads us to an
           important conclusion (obtained also in another way 
           in~\cite{VY95}) that the law of the gamma process is a limit 
           of renormalized stable processes when $\alpha\to0$ 
           (see Section~\ref{s:limit}).          
\end{description}       
           The last property allows us to consider the gamma process as an 
           antithesis to the Brownian motion: both correspond to the extreme values
           of the parameter $\alpha\in[0,2]$, 
           and in between we have 
           stable processes which can be considered as a deformation from the
           Brownian motion to the gamma process. 
           Moreover, we can see many strikingly similar properties,
           for example:
\begin{itemize}         
\item[1'.] Classical Bernstein's characteristic property of Gaussian 
         measures (components of both pairs of variables $(X,Y)$ and 
         $(X+Y,X-Y)$ are 
         independent if and only if the distributions of $X$ and $Y$ are 
         Gaussian) is similar to our property~1.
\item[2'.] Let $(B_t,\,t\ge0)$ be a Brownian motion.  
         For every $t>0$ the Brownian bridge $(B_u-\frac{u}{t}B_t,\, u \leq t)$
         is independent on $(B_v,\, v \geq t)$, which is similar to our property~2.
\item[3'.] For every  $b \in\Bbb R$ and every $t>0$ the law (a measure in the space of
         realizations) of the process $(B_u+bu,\, u \leq t)$ is equivalent to that of 
         $(B_u,\, u\leq t)$ with Radon--Nikodym density $\exp(bB_t-b^2 t/2)$.
         The large group of symmetries of the gamma process is a multiplicative 
         group of hyperbolic rotations (see Section~\ref{s:main}). 
         The law of the Brownian motion has a large orthogonal group of symmetries. 
         Stable processes perhaps also have such large groups of symmetries, but they
         consist of non-linear transformations, and can be 
         considered as a homotopy between hyperbolic and orthogonal rotations.
\end{itemize}                 
         
         In view of this parallelism, it would be quite interesting
         to connect better the Brownian motion and the gamma process, 
         and one such  attempt is recorded in the Appendix, where 
         we have gathered further properties of the gamma process which
         are not so close to the main context of the paper.

\begin{description}
\item[5. Generalized stability property:]
        In view of the previous property stating that the gamma process 
        is a limit of
        {\it stable} processes it is natural
        to ask what kind of stability do the gamma distribution
        and the gamma process have? In Section~\ref{s:stabletheory}
        we introduce a notion of generalized
        stability for {\it sigma-finite measures}, and show 
        that the Lebesgue measure on the real line and the quasi-Lebesgue
        measure which is a sigma-finite measure equivalent
        to the law of the gamma process are also stable in this new sense.
\end{description}
         
         Now we present our general framework and the list of topics 
         touched on in the paper. 
         The first one is the theory of Poisson--Dirichlet measures
         from the point of view of the gamma process. 
More general, we consider 
the theory of {\it L\'evy processes}
(a class of processes with independent
values) without Gaussian component. The laws of such processes
admit a canonical decomposition into a so-called {\it conic part},
i.e.~a measure on the cone of positive series, and a standard product
measure on sequences of points of the base space. This decomposition is of
the most general character. It appeared first in~\cite{FK72}, but the
original proof was rather complicated and shaded
the key relations with measures on positive
series and a special role of stable and gamma processes, while our proof is based
on very general and simple considerations. 

This decomposition leads us to
the theory of measures on 
positive series, i.e.~the theory of
{\it random positive series}, and their projections on the simplex of series
with unit sum. 
Particularly, we deal with the 
{\it Poisson--Dirichlet measures} $\operatorname{PD}(\th)$ 
and their generalizations which were studied by
Kingman~\cite{Ki75}, 
Vershik--Shmidt~\cite{VeSh77}, 
Pitman--Yor~\cite{PY97} and their followers, and have 
numerous applications in combinatorics (e.g.~\cite{VeSh77, Ha94}), 
number theory~\cite{Bi72, Ve86}, 
mathematical biology (e.g.~\cite{Ki78}, see a detailed survey in~\cite{Ew90}), etc.
See also some interesting discussion on the Poisson--Dirichlet 
         measures in most recently published  book~\cite{Bi99}.
These random series may be roughly characterized as ``the most random convergent
series'' --- a kind of white noise on convergent series.

The key feature in this work is the quasi-invariance of the
gamma measure (the law of the gamma process)
with respect to an infinite-dimensional group of
multiplicators. It was first discovered and used in the works of
Gelfand--Graev--Vershik~\cite{GGV73, GGV83, GGV85} 
on the representation theory of current groups, more
exactly, of the group $\operatorname{SL}(2,F)$, where $F$ is an algebra of functions on a
manifold. This property of the gamma measure followed from rather indirect
considerations. The same considerations prompted the existence of an
equivalent $\si$-finite measure which is invariant under multiplications by
non-negative functions with zero integral of logarithm. We call this measure
{\it quasi-Lebesgue}, because of its key property which is an infinite-dimensional
generalization of the well-known property 
of the finite-dimensional Lebesgue measure. We mean 
invariance of the Lebesgue measure
under the action of diagonal matrices
with determinant~$1$ (a Cartan subgroup). 
In this work we prove this quasi-invariance directly,
starting from the characteristic functional (the Laplace transform) of the
gamma measure, and construct explicitly the corresponding $\si$-finite
measure. The same quasi-invariance implies new symmetry properties of the
Poisson--Dirichlet measures 
with respect to Markovian transformations. 

However, the most important is the following link outlined in~\cite{VY95}. Both the
gamma measure and the quasi-Lebesgue measure are weak limits of
$\al$-stable measures (the laws
of $\al$-stable processes) when $\alpha$ goes to zero. In terms of the
Poisson--Dirichlet measures this fact was proved in~\cite{PY97}. More
exactly, Pitman--Yor~\cite{PY97} define a two-parameter family of
measures $\operatorname{PD}(\alpha,\theta)$ on the simplex of positive series with unit sum
and show that they
converge to the Poisson--Dirichlet measures $\operatorname{PD}(\th)$ when $\alpha\to0$. It turns out
that this convergence follows from the convergence of renormalized
$\alpha$-stable measures to the gamma measure when $\al\to0$. This
convergence is a ``commutative'' analogue of the key fact discovered
in~\cite{GGV73} which deals with the limit (more exactly,
the derivative with respect to the parameter in $0$ which corresponds to a
renormalization before taking the limit) 
of positive definite spherical functions of the
complementary series of $\operatorname{SL}(2,R)$ when the parameter tends to a critical value
(the so-called canonical state). Thus the Poisson--Dirichlet measures are
directly related to the representation theory.  

Another corollary of the quasi-invariance of the gamma measure allows to
obtain easily the {\it Markov--Krein identity}
which in our context
relates the distribution of a linear functional on the gamma
process and the distribution of the same functional on the normalized gamma
process.

\medskip
The paper is organized as follows.

In Section~\ref{s:levy} we define a general class of L\'evy processes. The
main properties of these processes are studied in Section~\ref{s:decomp}. 
Though in the sequel we deal only with the stable and gamma
processes, the Decomposition Theorem~\ref{th:decomp}
is proved by so general and natural considerations that we
give it in the most general form that does not complicate the argument. This
theorem states that
the law of a L\'evy process is the product of the conic part (the measure 
on the cone of positive convergent series) and a product measure on sequences 
of points of the base space. In
fact, our proof applies to even more general processes.
Theorem~\ref{th:conic} is a characterization of 
measures on the cone of positive convergent series
which are obtained as conic parts of L\'evy processes
(so-called measures of {\it product type}). In some sense, these measures 
enjoy the greatest possible
independence of coordinates. 
Note that passing from the simplex to the cone
(a kind of poissonization) simplifies many questions. For
example, characterization of the conic parts of L\'evy processes
is simpler than 
characterization of the simplicial parts. 

Section~\ref{s:gamma} contains definition and basic properties of the gamma process and
the Poisson--Dirichlet measures 
which are the simplicial parts of the gamma process. 

In Section~\ref{s:main} we prove the key property of the gamma
measure, namely its quasi-invariance with respect to a group of multiplicators
(Theorem~\ref{th:main}).
This group is rather wide and consists of all non-negative measurable
functions with summable logarithm. Thus the gamma measure is a
multiplicative analogue of the Wiener measure which is quasi-invariant with
respect to a wide group of additive shifts. This parallel should 
certainly be considered more thoroughly. It is not known  whether
the measure with such supply of preserving transformations is unique.
We give a partial counter-example for a smaller group of transformations. 
As to the quasi-invariance under the
whole group, the gamma measure seems to be the unique measure enjoying this
property.

In Section~\ref{s:PD} we apply
quasi-invariance of the gamma process to obtain the corresponding
property of the Poisson--Dirichlet measures.

In Section~\ref{s:sigma} we introduce a $\si$-finite measure which is
already invariant (projective invariant) under the same group of
multiplicators. It is called {\it quasi-Lebesgue} since it generalizes a
well-known property of the Lebesgue measure in $\Bbb R^n$. 

Section~\ref{s:stable} is devoted to definitions and basic properties of the
stable processes. The simplicial parts of these properties are related to the
two-parameter Poisson--Dirichlet measures $\operatorname{PD}(\al,\th)$.

In Section~\ref{s:limit} we prove the statement suggested in~\cite{VY95} 
that the gamma process is a weak limit of
renormalized stable processes. As was noted above, this fact is related to
the representation theory~\cite{GGV73, GGV83, GGV85}.

In Section~\ref{s:stabletheory} we give a new definition of stability which can
be applied to $\si$-finite measures as well (unlike the classical
definition). According to this definition, the quasi-Lebesgue measure is
zero-stable. 

Finally, in Sections~\ref{s:mk} and~\ref{s:genmk} we deal with the
{\it Markov--Krein transform} and its generalizations. We present a new
probabilistic interpretation of this transform: formulae of this kind relate
the distribution of a linear functional on the process with the distribution
of the same functional on the normalized process. This interpretation sets
the same question for general L\'evy processes.

In the Appendix we try to trace some further connection of our 
          topics, namely, some new and unexpected links with 
          Brownian motion.

The topics touched upon in this paper stimulate 
many new problems, only a small part of which is
mentioned above. 

\section{Definition of L\'evy processes on arbitrary spaces}\label{s:levy}  
Let $(X, \nu)$ be a standard Borel space with a non-atomic finite  
non-negative measure $\nu$, and let $\nu(X)=\theta$ be the total charge of $\nu$. 
We denote by 
$$
D=\left\{\sum z_i\de_{x_i},\;x_i\in X,\,z_i\in\Bbb R,\sum|z_i|<\infty \right\}
$$ 
a real linear space of all finite real discrete measures 
on $X$. 

Consider a class of measures $\La$ on the half-line $\Bbb R_+$ satisfying the 
following conditions,
\begin{eqnarray}
&&\La(0,\infty)=\infty, \label{cond1} \\
&&\La(1,\infty)<\infty, \label{cond2} \\
&&\int_0^1sd\La(s)<\infty, \label{cond3} \\
&&\La(\{0\})=0. \label{cond4}
\end{eqnarray}

Let $F_\La$ be the infinitely divisible
distribution with L\'evy measure $\La$, 
i.e.~the Laplace transform $\psi_\La$ of $F_\La$ is given by
$$
\psi_\La(t)=\exp\left(-\int_0^\infty(1-e^{-ts})d\La(s)\right).
$$ 

\begin{definition}
A {\it homogeneous L\'evy process} on the space $(X,\nu)$ with
L\'evy measure $\La$ satisfying~(\ref{cond1})--(\ref{cond4}) 
is a generalized process on $D$ whose law $P_\La=P_\La(\nu)$ has
the Laplace transform
\begin{equation}
\Bbb E\left[\exp\left(-\int_Xa(x)d\eta(x)\right)\right]=
\exp\left(\int_X\log\psi_\La(a(x))d\nu(x)\right),
\label{laplace}
\end{equation}
where $a$ is an arbitrary non-negative bounded Borel function on $X$.
\label{def:main}
\end{definition}

The correctedness of this definition is guaranteed by the following explicit construction
(see \cite[chapter~8]{Ki93}). 
Consider a Poisson point process on the space
$X\times\Bbb R_+$ with mean measure $\nu\times\La$.
We associate with a realization $\Pi=\{(X_i,Z_i)\}$ of this process an element
\begin{equation}
\eta=\sum_{(X_i,Z_i)\in\Pi} Z_i\de_{X_i}\in D. 
\label{poisson}
\end{equation}
Then $\eta$ is a random discrete measure obeying
the law $P_\La$. Note that if $\La$ is a $\de$-measure $\de_z$ for some $z\in\Bbb R_+$,
then $\Pi$ is the Poisson process on the set $X\times\{z\}$ (which we 
identify with $X$) with mean measure $\nu$, and the corresponding random element
$\eta$ is a measure that has equal charges $z$ at the points of this process.
Thus a L\'evy process with an arbitrary measure $\La$ is  
a continual convolution of independent Poisson processes on $X$ 
corresponding to different levels (charges).

It follows that the law $P_\La$ of the L\'evy process is concentrated on the cone 
$D^+=\{\sum z_i\de_{x_i}\in D:\;z_i>0\}\subset D$ consisting
of all finite positive discrete measures on X. The conditions~(\ref{cond1})--(\ref{cond4})
imposed on the measure $\La$
have the following meaning:~(\ref{cond1}) implies that the random
measure $\eta$ 
has an infinite number of atoms;~(\ref{cond2}) together with~(\ref{cond3})
guarantees that $\eta$ is a finite measure,
i.e.~the sum of charges is finite; 
finally,~(\ref{cond4}) means that our L\'evy process has no
Gaussian component.

{\bf Remarks. 1.}
Our definition of the L\'evy process is 
closely related to the notion of completely random measure, see~\cite{Ki67},
\cite[chapter~8]{Ki93}.

{\bf 2.}
If $X=\Bbb R_+$ and $\nu$ is the Lebesgue measure on $\Bbb R_+$, we recover an ordinary
definition of a {\it subordinator} (a process with stationary positive independent 
increments) corresponding to the L\'evy measure $\La$.

{\bf 3.}
It is easy to see that $P_\La(\nu)=P_{\th\La}(\nu/\th)$, i.e.~we may consider
only normalized parameter measures $\nu$. Thus {\it in the sequel we assume $\nu(X)=1$}.

\section{Decomposition theorem for L\'evy processes and
measures of product type on the cone}\label{s:decomp}
Consider the cone 
$$
C=\{z=(z_1,z_2,\ldots):\;z_1\ge z_2\ge\ldots\ge0,\,\sum z_i<\infty\}\subset l^1.
$$ 
We now define a special class of measures on $C$ indexed by
infinitely divisible distributions on the half-line.
Fix an integer $n\in\Bbb N$ and a probability vector
$p=(p_1,\ldots, p_n)$ 
(i.e.~a vector $p$ with $p_1,\ldots,p_n>0$ and $p_1+\ldots+p_n=1$).
Consider a sequence $\xi_i$ of i.i.d.~variables such that 
$P(\xi_i=k)=p_k$ for $k=1,\ldots n$. For $Q=(Q_1,Q_2,\ldots)\in C$, denote by 
$\Si_k^{(p)}=\Si_k^{(p)}(Q)$
the random sum $\Si_k^{(p)}=\sum_{i:\xi_i=k}Q_i$. 
Let $Q$ be a random series with distribution 
$\varkappa$ on $C$ such that 
the distribution $F$ of the sum $\sum Q_i$ is infinitely 
divisible.

\begin{definition}  
We say that a series $Q$ (and its distribution
$\varkappa$) is of {\it product type}, if for each 
$n\in\Bbb N$ and each probability vector $p$
the sums $\Si_1^{(p)},\ldots\Si_n^{(p)}$ are independent and $\Si_k^{(p)}$ obeys the law
$F^{*p_k}$. 
\end{definition}

We define a map $T:D^+\to C\times X^\infty$ by
\begin{equation}
T\eta=\big((Q_1,Q_2,\ldots),\;(X_1,X_2,\ldots)\big),\quad\text{ if }\quad
\eta=\sum Q_i\de_{X_i}.
\label{T}
\end{equation}

\begin{definition}
Let $P$ be a distribution on the space $D^+$, and 
let $\eta$ be a random process obeying the law $P$. 
The random sequence of charges $Q_1,Q_2,\ldots$ is called the
{\it conic part} of the process $\eta$, and its distribution  
on the cone $C$ is called
the {\it conic part} of the law $P$.
\end{definition} 

Note that in view of representation~(\ref{poisson}) 
the conic part of the L\'evy process
with L\'evy measure $\La$ is just the ordered sequence of points of the Poisson process
on $\Bbb R_+$ with mean measure $\La$. Thus the conic part depends only on $\La$ and not 
on the $(X,\nu)$. In fact, the following theorem shows that studying the L\'evy process 
may be essentially reduced   
to studying its conic part, since the construction of the process includes
the parameter measure in a trivial way. 
This fundamental property of homogeneous L\'evy processes
is a particular case of the representation theorem first
proved in \cite{FK72}. We present here a simpler proof of this fact. 
The key point of our proof is the following lemma.
Let $(X,\nu)$ be a standard Borel space with continuous probability measure $\nu$.
Denote $X^k=X\times\ldots\times X$ ($k$ factors), 
$\nu^k=\nu\times\ldots\times\nu$ ($k$ factors) and let
$\nu_{diag}$ be the image of $\nu$ under the diagonal map $x\to(x,\ldots,x)$. 

\begin{lemma}\label{l:inv} 
Let $\tau$ be some continuous probability measure on $X^k$.
If for each measure preserving 
transformation $L$ of $(X,\nu)$, the transformation $L^k=L\times\ldots\times L$
($k$ factors) preserves $\tau$, 
then $\tau$ is a convex combination of $\nu^k$ and $\nu_{diag}$. 
\end{lemma}

\begin{proof}
For simplicity, assume $k=2$, the general case being quite similar.
The diagonal $\De=\{(x,x),\,x\in X\}$ is obviously an invariant subset
for the group $G=\{L\times L\}$, where $L$
runs over the set of all $\nu$-preserving transformations of the space $X$. Thus
it suffices to show that if $\tau$ 
is concentrated on the set $(X\times X)\setminus \De$, then
$\tau=\nu\times\nu$, and if $\tau$ is concentrated on $\De$, then $\tau=\nu_{diag}$.

In the first case let $\xi_n=\{A_i\}_{i=1}^{2^n}$ be an arbitrary partition 
of the space $X$ into $2^n$ sets of equal $\nu$-measure $1/2^n$. Denote by $\xi_n^2$
the corresponding partition of the space $X\times X$, 
i.e.~$\ti\xi_n=\{A_{ij}\}$, where $A_{ij}=A_i\times A_j$. 
The group $G$ acts 
transitively on the set of non-diagonal elements of $\xi_n$.
Thus all non-diagonal
elements have equal $\tau$-measure. Denote 
$Y_n=(X\times X)\setminus\cup(A_i\times A_i)$ and $\eps_n=\tau(Y_n)$. 
Since $\tau$ is concentrated on $(X\times X)\setminus \De$, 
we have $\eps_n\to1$ as $n\to\infty$. Considering finer partitions and using
the above argument, we obtain that for each $k\in\Bbb N$,  
if a rectangle $A\times B\subset Y_n$ and $\nu(A)=\nu(B)=\nu(Y_n)/2^k$, then
$\tau(A\times B)=\eps_n/4^k$. But then the restriction of $\tau$ on the set $Y_n$ equals
$\eps_n\cdot(\nu\times\nu)$. Letting $n\to\infty$, we obtain $\tau=\nu\times\nu$.

In the second case, identifying the diagonal $\De$ with $X$, we obtain that $\tau$ is a 
measure on $X$ which is invariant under all $\nu$-preserving transformations, 
hence obviously $\tau=\nu$.
\end{proof}

\begin{theorem}\label{th:decomp} 
Let $\eta=\sum Q_i\de_{X_i}$ be a homogeneous L\'evy process on the space 
$(X,\nu)$ with L\'evy measure $\La$. 
Then $TP_\La=\varkappa_\La\times\nu^\infty$, 
i.e.~$X_1,X_2,\ldots$ is a sequence of i.i.d.~random variables
with common distribution $\nu$, and this sequence is independent of 
the conic part
$\{Q_i\}_{i\in\Bbb N}$.  
\end{theorem}      

\begin{proof}
Let $L:X\to X$ be a $\nu$-preserving transformation of $X$. This transformation
acts on the space $D$ by substituting coordinates, 
i.e.~$\sum z_i\de_{x_i}\mapsto\sum z_i\de_{Lx_i}$, and 
it is clear that the law $P_\La$ of a  
homogeneous L\'evy process on $(X,\nu)$ is invariant under $L$.
Denote by $P_\La^z$ 
the conditional measure of $P_\La$ given the conic part equal to $z\in C$.
The transformation $L$ acts ``fibre-wise'', i.e.~it does not
change the conic part, hence $L$ preserves almost all
conditional measures $P_\La^z$. In particular, if we denote by $(P_\La^z)_k$
the conditional distribution of the first $k$ points $X_1,\ldots,X_k$
on the space $X^k$, then
the transformation $L^k$ preserves $(P_\La^z)_k$.

Now it follows from Lemma~\ref{l:inv} that for almost all $z$ 
$(P_\La^z)_k=\nu^k$ for all $k$, i.e.~$P_\La^z=\nu^\infty$,
and Theorem  1 follows.
\end{proof}
 
\begin{theorem}\label{th:conic}
The measure $\varkappa$ on the cone $C$ is the conic part of some L\'evy process $P_\La$ 
with L\'evy measure $\La$ satisfying~(\ref{cond1})--(\ref{cond4})
if and only if it is of product type with $F=F_\La$. 
\end{theorem}

\begin{proof}
Given fixed 
$n\in\Bbb N$ and a probability vector $p=(p_1,\ldots,p_n)$, consider a partition 
$X=A_1\cup\ldots\cup A_n$ of the space $X$ such that $\nu(A_k)=p_k$, $k=1,\ldots,n$. 
Let $X_1,X_2,\ldots$ be the sequence of i.i.d.~variables with common distribution
$\nu$ and assume $\xi_i=k$, if $X_i\in A_k$.
Then the random variables $\xi_i$ form a sequence of i.i.d.~variables, and
$P(\xi_i=k)=p_k$. 
Consider a random process 
\begin{equation}
\eta=\sum Q_i\de_{X_i}, 
\label{eta}
\end{equation}
where the sequence $Q_1,Q_2,\ldots$ is independent of 
$\{X_i\}$ and obeys the law $\kappa$. Let $\Si_k^{(p)}=\sum_{i:\xi_i=k}Q_i$.
It is easy to see that
for arbitrary $t_1,\ldots,t_n>0$ 
\begin{equation}
\Bbb E\left[\exp\left(-\sum_{k=1}^n t_k\Si_k^{(p)}\right)\right]=
\Bbb E\left[\exp\left(-\int_Xa(x)d\eta(x)\right)\right],
\label{conn}
\end{equation}
where $a$ is a step function such that $a(x)=t_i$, if $x\in A_i$. 

Now let $\kappa$ be the conic part of the law $P_\La$ of some L\'evy process.
Then, by Theorem~\ref{th:decomp}, the process $\eta$ defined by~(\ref{eta})
obeys $P_\La$, and
it follows 
from the Laplace transform formula~(\ref{laplace}) that the right-hand side 
of~(\ref{conn}) equals
$$
\prod_{i=1}^n \psi_\La(t_i)^{\nu(A_i)}=\prod_{i=1}^n \psi_\La(t_i)^{p_i}.
$$ 
Since the left-hand side of~(\ref{conn}) is the Laplace transform of the common
distribution of the variables $\Si_1^{(p)},\ldots,\Si_n^{(p)}$, we obtain 
that $\Si_1^{(p)},\ldots,\Si_n^{(p)}$ are independent, and $\Si_k^{(p)}$
obeys the law $F^{*p_k}$, i.e.~$\kappa$ is of product type.

Conversely, let $\kappa$ be a measure of product type corresponding to
an infinitely divisible law $F$. Define a random process $\eta$ 
on an arbitrary measurable space $(X,\nu)$ satisfying the conditions of
Definition~\ref{def:main} by~(\ref{eta}).
The above argument shows that $\eta$ satisfies~(\ref{laplace})
with $\La$ equal to the L\'evy measure of $F$
for all positive step functions $a$, and one can easily extend this to all bounded
positive Borel functions by continuity. Thus $\eta$ obeys $P_\La$, and $\kappa$
is the conic part of $P_\La$.
\end{proof}

The strong law of large numbers for Poisson processes (see \cite[4.55]{Ki93})
implies the following 
result on the asymptotic behavior of the vectors $Z=(Z_1,Z_2,\ldots)\in C$
obeying the law $\kappa_\La$.

\begin{proposition}[\cite{Ki75}]\label{p:asymp}
Let $m_\La(t)=\La(t,\infty)$, $t>0$. Then
\begin{equation}
\lim_{n\to\infty}\frac{m_\La(Z_n)}{n}=1
\label{asymp}
\end{equation} 
for almost all with respect to $\kappa_\La$ vectors $Z\in C$.
\end{proposition}

Note that we may rewrite~(\ref{asymp}) as 
\begin{equation}
\lim_{z\to0}\frac{m_\La(z)}{\#\{i:\;Z_i>y\}}=1.
\label{asymp'}
\end{equation}
In other words, if $\eta$ is a random discrete measure obeying
the law $P_\La$, then the number of charges of $\eta$ which are greater than $z$ has
the same asymptotics when $z\to0$ as the tail $\La(z,\infty)$ of the L\'evy measure $\La$.  

Denote by $D^+_1\subset D^+$ the simplex 
of all normalized atomic measures. Then 
$D^+ =D^+ _1 \times [0, \infty)$, i.e.~each $\eta \in D^+$ can be represented as
\begin{equation}
\eta=({\eta}/{\eta(X)}, \eta(X)). 
\label{normaliz}
\end{equation}
The second coordinate 
is the total charge of the measure $\eta$. It follows from the definition 
of a L\'evy process that $\eta(X)$ obeys the infinite divisible law $F_\La$ 
corresponding to the L\'evy measure $\La$.
The first coordinate 
is called the {\it normalization} of the measure $\eta$. 
Note that, in general, the law of a L\'evy process is not a product measure
in this decomposition (see Lemmas~\ref{l:ind} and~\ref{l:Lu} below). 

Using this decomposition,
consider a map
$T':D^+\to\Bbb R_+\times\Si\times X^\infty$, where 
$
\Si=\{y=(y_1,y_2,\ldots):\;y_1\ge y_2\ge\ldots\ge0,\; y_1+y_2+\ldots=1\}
$
is the infinite-dimensional simplex, and
$$
T'\eta=\big(\eta(X),\;(Q_1/\eta(X),Q_2/\eta(X),\ldots),\;
(X_1,X_2,\ldots)\big),\quad\text{if}\quad\eta=\sum Q_i\de_{X_i}.
$$

\begin{definition}
The normalized sequence of charges $Q_1/\eta(X),Q_2/\eta(X),\ldots$
is called the {\it simplicial part} of the process and its distribution $\si_\La$
on $\Si$ is called the {\it simplicial part} of the law $P_\La$.
\end{definition}

\section{The gamma process and Poisson--Dirichlet distributions}\label{s:gamma}

Let $(X, \nu)$ be a standard Borel space with a non-atomic finite  
non-negative measure $\nu$, and let $\nu(X)=\theta$ be the total charge of $\nu$.

\begin{definition} 
{\it The gamma process} with scale parameter $\la>0$ on the space $(X,\nu)$ 
is a L\'evy process on $(X,\nu)$ corresponding to the L\'evy measure 
with density $d\La^\la_\Ga(z)=z^{-1}e^{-\la z}dz$, $z>0$.
\end{definition}

As in general case, 
instead of considering a non-normalized parameter measure $\nu$,
we may take its normalization $\bar\nu=\nu/\nu(X)$ and the L\'evy measure with
density $\th z^{-1}e^{-\la z}dz$. 
But in this particular case it is often more convenient to use the above definition
with a non-normalized measure. In the sequel we shall use both variants without
additional mention. 

The corresponding infinitely divisible law is
the gamma distribution ${\cal G}_{\th,\la}$
on $\Bbb R_+$ with shape parameter $\th$ and scale
parameter $\la$, i.e. 
$$
d{\cal G}_{\th,\la}=\frac{\la^\th}{\Ga(\th)}t^{\th-1}e^{-\la t}dt,\quad t>0.
$$

Note that $\la$ is a trivial scale parameter. Namely, if $\eta$ is the gamma process  
with scale parameter $1$, then the gamma process $\eta^\la$
with scale parameter $\la$ is obtained from $\eta$ by multiplying by $\la$, 
i.e.~$\eta^\la=\la\eta$. Thus {\it we will consider only gamma processes with scale 
parameter~$1$}. 

The law $P_\Ga=P_\Ga(\nu)$ of the gamma process
(called the {\it gamma measure}
on the space $(X,\nu)$) is thus
given by the Laplace transform 
\begin{equation}
\Bbb E_{\Ga}\left[\exp\left(-\int_X a(x)d\eta(x)\right)\right]=
\exp\left(-\int_X\log\left(1+a(x)\right)d\nu(x)\right),
\label{galapl}
\end{equation}
where $a$ is an arbitrary non-negative bounded Borel function on $X$. 
 
Let ${\cal M}={\cal M}(X,\nu)$ be the set of (classes\!$\mod 0$ of)
non-negative measurable functions on the space $X$
with $\nu$-summable logarithm,
$$
{\cal M}=\left\{a:X\to\Bbb R_+:\;\int_X|\log a(x)|d\nu(x)<\infty\right\}.
$$
It follows from the Poisson construction~(\ref{poisson})  
and Campbell's theorem on sums over Poisson processes (see~\cite[3.2]{Ki93})
that each function $a\in{\cal M}$ correctly defines a measurable linear functional
$\eta\mapsto f_a(\eta)=\int_X a(x)d\eta(x)$ on $D$, and formula~(\ref{galapl}) holds for
all $a\in{\cal M}$.

It is well known that the gamma distribution enjoys the following property. If
$Y$ and $Z$ are independent gamma
variables with the same scale parameter, 
then the variables $Y+Z$ and $\frac Y{Y+Z}$ are independent. Moreover, a remarkable
result of Lukacs~\cite{Lu65} (similar to the famous Bernstein's characterization
of normal distributions) 
states that this property is characteristic of the
gamma distribution, i.e.~if $Y$ and $Z$ are independent non-degenerate positive
random variables, and the variables $Y+Z$ and $\frac Y{Y+Z}$ are independent, then
both $Y$ and $Z$ have gamma distributions with the same scale parameter.
In other words, the independence condition may be formulated as follows. Let us describe
a point $x=(x_1,x_2)$ in the first quadrant $\Bbb R_+\times\Bbb R_+$ by the sum
$x_1+x_2$ of its coordinates and its projection onto the unit 2-simplex (i.e.~interval)
$\{y=(y_1,y_2):\,y_1,y_2\ge0,\,y_1+y_2=1\}$. 
Then the distribution ${\cal G}_{\th_1,\la}\times{\cal G}_{\th_2,\la}$ 
on $\Bbb R_+\times\Bbb R_+$
is a product measure in these coordinates.
These results imply the corresponding statements for the gamma process which are a key
point for many important properties of $P_\Ga$.

\begin{lemma} 
\label{l:ind}
In representation~(\ref{normaliz}) the gamma measure is a product measure
$P_\Gamma = {\cal G}_\th\times\bar P_\Gamma$, i.e.~the total charge $\ga(X)$
of the gamma process and the normalized gamma process $\bar\ga=\ga/\ga(X)$ are independent.
The distribution of the total charge is the gamma distribution ${\cal G}_{\th,1}$
with shape parameter $\th$ and scale
parameter 1.
\end{lemma}

\begin{lemma}
\label{l:Lu} 
If the law $P_\La$ of some 
L\'evy process is a product measure in representation~(\ref{normaliz}), then
$P_\La$ is a gamma process, i.e.~$d\La(z)=z^{-1}e^{-\la z}dz$, $z>0$,
for some $\la>0$.
\end{lemma}
    
\begin{definition}[\cite{Ki75}]\label{def:PD}
The simplicial part of the gamma 
measure $P_\Ga(\nu)$ with $\nu(X)=\th$ is called the {\it Poisson--Dirichlet}
distribution with parameter $\th$ and denoted by $\operatorname{PD}(\th)$.
\end{definition}

The above definition is just one of many other possible definitions of 
the Poisson--Dirichlet distributions. These distributions arise in many fields
of pure and applied mathematics. They play an important role in statistics because
of their connection with Dirichlet distributions~\cite{Ki75} and Dirichlet random 
measures~\cite{Fe73}. In number theory Poisson--Dirichlet measures arise
in the problem of distribution of prime divisors of a random integer~\cite{Bi72,Ve86}.
There are many asymptotic
combinatorial problems leading to the measures $\operatorname{PD}(\th)$, such as
the distribution of the cycle lengths of a random permutation~\cite{VeSh77} or
the distribution of the degrees of the 
irreducible factors of a random monic polynomial over a finite field~\cite{Ha94}, etc. The 
Poisson--Dirichlet distributions also play an important role in applications to
population genetics, ecology and physics. A (non-complete) survey of different aspects
of Poisson--Dirichlet measures can be found in~\cite{EwTa97}. 

It follows from Lemma~\ref{l:ind} that 
$T'P_\Ga={\cal G}_\th\times \operatorname{PD}(\th)\times\bar\nu^\infty$, 
i.e.~the conic part of the gamma measure is
a product measure ${\cal G}_{\th,1}\times \operatorname{PD}(\th)$. We call this measure 
the {\it conic Poisson--Dirichlet distribution} with parameter $\th$
and denote it by $\operatorname{CPD}(\th)$. 

Note that in case of the gamma process
$m(t)=\th\int_t^\infty s^{-1}e^{-s}ds\sim -\th\log t$. Thus,
by Proposition~\ref{p:asymp},
\begin{equation}
\lim_{n\to\infty}\frac{\log Z_n}{n}=-\frac{1}{\th}
\label{PDasymp}
\end{equation}
almost surely with respect to $\operatorname{CPD}(\th)$. It follows that the same asymptotics holds
for $\operatorname{PD}(\th)$, i.e.~$\lim_{n\to\infty}\frac{\log Y_n}{n}=-\frac{1}{\th}$ for almost
all  vectors $Y\in\Si$ with respect to $\operatorname{PD}(\th)$. In particular, we see that 
the measures $\operatorname{PD}(\th)$ (as well as $\operatorname{CPD}(\th)$) are mutually singular for different $\th$.

Many properties of ordinary Poisson--Dirichlet distributions
have their natural analogues for conic distributions. For example, it is well known
that the measure $\operatorname{PD}(\th)$ may be obtained by the following {\it stick breaking process}.
Let $Y_1$ be a random variable on the interval $[0,1]$ obeying the law
$\th(1-t)^{\th-1}dt$, $t\in[0,1]$. If we have already constructed variables 
$Y_1,\ldots,Y_n$, then $Y_{n+1}$ has the same distribution scaled on the 
interval $[Y_n,1]$. Thus we obtain a random sequence $0=Y_0<Y_1<Y_2<\ldots<1$.
Let $Z_k=Y_k-Y_{k-1}$, $k=1,2,\ldots$.
The Poisson--Dirichlet measure $\operatorname{PD}(\th)$ is the distribution
of the order statistics $Z_{(1)}\ge Z_{(2)}\ge\ldots$ of the sequence 
$Z_1,Z_2,\ldots$. It follows that the conic Poisson--Dirichlet measure $\operatorname{CPD}(\th)$ may
be obtained by the randomized version of this procedure. Namely, for the first step we
choose the random length $L$ of the interval with gamma distribution ${\cal G}_{\th,1}$,
and then proceed as before starting with the interval $[0,L]$.  

As shown in~\cite{PY97}, the Poisson--Dirichlet measures $\operatorname{PD}(\th)$ can be
naturally included into a two-parameter family $\operatorname{PD}(\al,\th)$ of distributions on the 
simplex $\Si$. This family is obtained by the following
non-stationary version of the stick breaking process.
Let $Y_1$ be a random variable on the interval $[0,1]$ obeying the beta distribution
$B(1-\al,\th+\al)$. If we have already constructed variables 
$Y_1,\ldots,Y_n$, then $Y_{n+1}$ has the beta distribution $B(1-\al,\th+(n+1)\al)$
scaled on the 
interval $[Y_n,1]$. Thus we obtain a random sequence $0=Y_0<Y_1<Y_2<\ldots<1$.
Let $Z_k=Y_k-Y_{k-1}$, $k=1,2,\ldots$.
The two-parameter Poisson--Dirichlet measure $\operatorname{PD}(\al,\th)$ is the distribution
of the order statistics $Z_{(1)}\ge Z_{(2)}\ge\ldots$ of the sequence 
$Z_1,Z_2,\ldots$.

The range of admissible parameters is the union of the sets 
$$
\{(\al,\th):\;0\le\al<1,\,\th>-\al\}\quad\text{and}\quad
\{(\al,-m\al):\;\al<0,m\in\Bbb N\}.
$$
The first case $\al\in(0,1)$ is the most interesting, the second one being
a sort of degenerate case. The ordinary Poisson--Dirichlet distributions $\operatorname{PD}(\th)$
correspond to $\al=0$, i.e.~$\operatorname{PD}(\th)=\operatorname{PD}(0,\th)$.
See~\cite{PY97} for various properties of the measures
$\operatorname{PD}(\al,\th)$. In particular, the distributions $\operatorname{PD}(\al,\th)$ with fixed $\al\ne0$ and
different $\th$ are absolutely continuous (unlike the case $\al=0$).
We discuss some questions related to the two-parameter Poisson--Dirichlet distributions
in Sections~\ref{s:stable} and~\ref{s:genmk}.

\section{Multiplicative quasi-invariance of the gam\-ma process}
\label{s:main}

As was mentioned above, the law $P_\La$ of each L\'evy process 
is invariant under all $\nu$-preserving 
transformations of the space $(X,\nu)$
which act on $D$ by substituting the coordinates.
However, the gamma measure $P_\Gamma$ enjoys additional
invariance properties.   
We present now a large group of linear transformations of the space $D$ 
(preserving the cone $D^+$) for which $P_\Gamma$ is a quasi-invariant measure. 
    
Consider the above defined class ${\cal M}$ of non-negative functions on $X$ 
with $\nu$-summable logarithm.
Each function $a\in{\cal M}$ defines not only a linear functional $f_a$ on $D$ but also
a multiplicator $M_a:D\to D$ by $(M_a\eta)(x)=a(x)\eta(x)$, that is 
$M_a\eta=\sum a(x_i)z_i\de_{x_i}$ for $\eta=\sum z_i\de_{x_i}$.
Note that the set $\cal M$ is a commutative group with respect to pointwise
multiplication of functions, and $M_a$ is a group action of $\cal M$.
Denote by $\ti a$ the function $\ti a(x)=(1/a(x))-1$.

The following property of the gamma process was first discovered 
in~\cite{GGV83, GGV85} in quite different terms; 
it plays an important role in the representation theory  
of the current group $\operatorname{SL}(2,F)$, where $F$ is the space of functions on
a manifold.
        
\begin{theorem}\label{th:main}
For each $a\in{\cal M}$, the measure $P_\Ga$ is
quasi-invariant under $M_a$, and the corresponding density is given by
the following formula,
\begin{equation}
\frac{d(M_aP_\Ga)}{dP_\Ga}(\eta)=\exp\left(-\int_X\log a(x)d\nu(x)\right)\cdot
\exp\left(-\int_X\ti a(x)d\eta(x)\right).
\label{density}
\end{equation}
\end{theorem}     

\begin{proof}
Fix $a\in{\cal M}$ and let $\xi=L_a\eta$. 
Consider an arbitrary function $b\in{\cal M}$. Then 
$\int_X b(x)d\xi(x)=\int_X a(x)b(x)d\eta(x)$. Thus, in view of~(\ref{galapl}),
the Laplace transform
$\Bbb E\left[\exp\left(-\int_X b(x)d\xi(x)\right)\right]$
equals
\begin{multline*}
\Bbb E\left[\exp\left(-\int_X a(x)b(x)d\eta(x)\right)\right]=
\exp\left(-\int_X\log\big(1+a(x)b(x)\big)d\nu(x)\right)=\\=
\exp\left(-\int_X\log a(x)d\nu(x)\right)\cdot
\exp\left(-\int_X\log\left(\frac1{a(x)}+b(x)\right)d\nu(x)\right).      
\end{multline*}
Using~(\ref{galapl}) once more, we may consider the last factor as the Laplace transform
of $P_\Ga$ calculated on the function $(\frac1{a(x)}-1)+b(x)=\ti a(x)+b(x)$.
Denote $I(a)=\int_X\log a(x)d\nu(x)$. Then we have
\begin{multline*}
\Bbb E\left[\exp\left(-\int_X b(x)d\xi(x)\right)\right]=I(a)\cdot
\Bbb E\left[\exp\left(-\int_X\left(\ti a(x)+b(x)\right)d\eta(x)\right)\right]=\\=
\Bbb E\left[I(a)\cdot\exp\left(-\int_X\ti a(x)d\eta(x)\right)\cdot
\exp\left(-\int_X b(x)d\eta(x)\right)\right],
\end{multline*}
and Theorem 3 follows.
\end{proof}   

In particular, if we consider multiplication by constant $c>0$, then the corresponding 
density depends only on the total charge, namely
\begin{equation}
\frac{d(M_cP_\Ga)}{dP_\Ga}(\eta)=\frac1{c^\th}\cdot
\exp\left(\left(1-\frac1c\right)\eta(X)\right).
\label{constdens}
\end{equation}

\begin{theorem}\label{th:erg}
The action of the group $\cal M$ 
on the space $(D^+,P_\Gamma)$  is ergodic.
\end{theorem}
\begin{proof}
Let $G:D^+\to\Bbb R$ be a $P_\Ga$-measurable functional on $D^+$
which is invariant under all $M_a$ i.e.~$G(M_a\eta)=G(\eta)$ a.e. with respect to 
$P_\Ga$. Consider an arbitrary Borel function  $k:\Bbb R\to\Bbb R$. Then
for each $a\in{\cal M}$
\begin{multline*}
\Bbb E\big[k(G(\eta))\big]=\Bbb E\big[k(G(M_a\eta))\big]=\\\Bbb E\left[k(G(\eta))\cdot 
\exp\left(-\int_X\ti a(x)d\eta(x)\right)\right]\cdot
\exp\left(-\int_X\log a(x)d\nu(x)\right),
\end{multline*}
where $\Bbb E$ denotes the expectation with respect to $P_\Ga$.
But in view of~(\ref{galapl}) the last factor equals
$$
\left(\Bbb E\left[\exp\left(-\int_X\ti a(x)d\eta(x)\right)\right]\right)^{-1}=
\left(\Bbb E\left[\exp\left(-f_{\ti a}(\eta)\right)\right]\right)^{-1},
$$
hence we have
$$
\Bbb E\left[k(G(\eta))\exp\left(-f_{\ti a}(\eta)\right)\right]=
\Bbb E\big[k(G(\eta))\big]\cdot
\Bbb E\left[\exp\left(-f_{\ti a}(\eta)\right)\right].
$$
Thus $G$ is independent of every functional $f_a$, and Theorem 4 follows.
\end{proof}

It is natural to ask whether the quasi-invariance
property stated in Theorem~\ref{th:main} is characteristic of the gamma process. 
If we fix the density, then the answer is positive (\cite{GGV85}), i.e.~the 
gamma measure is the only measure on $D^+$ 
satisfying~(\ref{density}). On the other hand,  
for a smaller group of multiplications, the answer may be negative,
as the following example shows.
Let us call a process {\it quasi-multiplicative} if its law is invariant 
under all transformations $M_a$ with {\it constant} $a>0$. For simplicity, 
consider {\it subordinators} (i.e.~L\'evy 
processes for $X=\Bbb R_+$ and $\nu$ equal to the Lebesgue measure).

\begin{proposition}[\cite{VY95}]\label{p:qm}
Let $\eta$ be a subordinator
with L\'evy measure $\La(dx)=k(x)dx$, where $k(x)>0$, and denote
$g(x)=xk(x)$. Then $\eta$
is quasi-multiplicative if and only if for all $a>0$
\begin{equation}
\int_0^1\left(\sqrt{g(x\,/\,a)}-\sqrt{g(x)}\right)^2\frac{dx}x<\infty.
\label{quasimult}
\end{equation}
\end{proposition}

It is shown in \cite{VY95} that for each $m<1/2$ any function $k_m(x)$ that
satisfies
\begin{equation}
k_m(x)=\frac1x\left(\log\frac1x\right)^{2m}\text{ for \ \ }x<1/2
\qquad\text{ and }\qquad
\int_{1/2}^\infty k_m(x)dx<\infty
\label{example}
\end{equation}
provides an example of a quasi-multiplicative subordinator that is not equivalent to any
scaled gamma process.
It is clear that if $\eta$ is quasi-multiplicative,
then its law is quasi-invariant under all step functions
with finitely many steps, but the quasi-invariance 
under the whole group ${\cal M}$ does not take place in this example.
It is not known if there exists a measure different from $P_\Gamma$
which has this property.

Computation of the law of hitting time of the drifted gamma process 
is coherent with the quasi-invariance of the gamma process. This fact is also
parallel to some property of Brownian motion and will be considered
in more details elsewhere.

\section{Quasi-invariance of the Poisson--Dirichlet distributions}\label{s:PD}

Let $a\in{\cal M}$.
According to the general theory of polymorphisms (see \cite{Ve77}),
the transformation $M_a$ induces a Markovian operator $R_{a}$ on the cone $C$.
Namely, let $z=(z_1,z_2,\ldots)\in C$. Consider the conditional
distribution $P^z_\Ga$ of the gamma process on $(X,\nu)$, 
given the conic part equal to
$z$. Then the random image of the point $z$ under
$R_{a}$ is the conic part
of the process $M_a\eta$, where $\eta$ 
obeys the law $P_\Ga^z$. It follows from Theorem~\ref{th:decomp} that 
$$
R_{a}z=V({a(X_1)z_1},{a(X_2)z_2},\ldots),
$$
where $(X_1,X_2,\ldots)$ is a sequence of i.i.d.~random variables on $X$
with common distribution $\nu$,
and  $V$ denotes a map that
arranges the coordinates in non-increasing order.

In a similar way, the transformation $M_a$ induces
a Markovian operator $S_{a}$ on the simplex $\Si$,
$$
S_{a}y=V\left(\frac{a(X_1)y_1}{\si},\frac{a(X_2)y_2}{\si},\ldots
\right),
$$
where the sequence $(X_1,X_2,\ldots)$ is as before, and 
$\si=a(X_1)y_1+a(X_2)y_2+\ldots$.

Note that the definitions of the operators $S_a$ and $R_a$ depend only on 
the distribution of the function $a$. Thus 
we may assume that $X=[0,1]$ and
$\nu=\th\la$, where $\la$ is the Lebesgue measure on the interval.

Theorems~\ref{th:main},~\ref{th:erg} immediately imply

\begin{theorem}\label{pdquasi}
{\rm 1)} The Poisson--Dirichlet distribution $\operatorname{PD}(\th)$
is quasi-invariant under the Markovian operator $S_{a}$
for all $a\in\cal M$, and 
$$
\frac{dS_{a}\operatorname{PD}(\th)}{d\operatorname{PD}(\th)}(y)=
\exp\left(-\th\int_0^1\log a(s)ds\right)\cdot
\int_0^\infty\frac{\si^{\th-1}}{\Ga(\th)}
\left(\prod_{i=1}^\infty L_{1/a}(\si y_i)\right)d\si,
$$
where $L_{1/a}(\cdot)$ is the Laplace transform of the distribution
of the function $1/a(t)$ with respect to the uniform
distribution on the interval $[0,1]$.

{\rm 2)} The Poisson--Dirichlet distribution $\operatorname{PD}(\th)$ is
ergodic with respect to $\{S_{a}\}_{a\in\cal M}$.
\end{theorem}

\section{Quasi-Lebesgue measure and a re\-pre\-sen\-ta\-ti\-on of the current group} 
\label{s:sigma}

In this section we define, following \cite{GGV83, GGV85},
a $\si$-finite measure on $D^+$ which is equivalent to $P_\Ga$ and 
invariant under a subgroup ${\cal M}_0$ of ${\cal M}$ consisting
of all functions $a\in{\cal M}$ such that
$\int_X\log a(x)d\nu(x)=0$. 

\begin{definition} 
Consider a $\si$-finite measure $\ti P_\Ga$ on $D^+$ defined by 
\begin{equation}
\frac{d\ti P_\Ga}{dP_\Ga}(\eta)=\exp(\eta(X)).
\label{sigmadens}
\end{equation}
It is called the quasi-Lebesgue measure. 
\end{definition}

Theorem~\ref{th:main} implies 

\begin{theorem}\label{th:sigma}
For each $a\in\cal M$, the quasi-Lebesgue measure $\ti P_\Ga$ is quasi-invariant under
$M_a$ with a constant density
$$
\frac{dM_a(\ti P_\Ga)}{d\ti P_\Ga}=\exp\left(-\int_X\log a(x)d\nu(x)\right).
$$
\end{theorem}

\begin{corollary}\label{c:inv}
If $\int_X \log a(x)d\nu(x)=0$, then the quasi-Lebesgue measure
$\ti P_\Gamma$ is invariant
with respect to $M_a$. 
\end{corollary}

We see that the measure  $\ti P_\Gamma$ is invariant with respect to an 
infinite-dimensional multiplicative group whose action 
generalizes
the action of the group of diagonal matrices with determinant $1$ in
a finite-dimensional vector space. Thus we may consider the measure 
$\ti P_\Gamma$ as an infinite-dimensional analogue of the Lebesgue measure.
This property was much used in 
\cite{GGV83, GGV85} for the representation theory 
of the group $\operatorname{SL}(2,F)$. 

Let us consider the group of triangular matrices of order~$2$ 
$$
T_{a,b}=\left(
\begin{array}{cc}
a(\cdot)& b(\cdot)\\
0&a(\cdot)^{-1}
\end{array}
\right)
$$
with $b,\log a\in L^1(X,m)$ for some
measurable space $(X,m)$ (a current group in the terminology of physicists).

\begin{theorem}
The formula
$$
U(T_{a,b})F(\eta)= 
\exp\left(\int_X\log a(x)d\nu(x)+i\int_X a(x)b(x)d\eta(x)\right)F(M_a\eta).
$$ 
defines a unitary irreducible representation of this group
        in the space $L^2(\ti P_\Gamma)$.
\end{theorem}

\begin{proof}
The representation is correctly defined and its unitarity follows from
the invariance property of $\ti P_\Ga$. The irreducibility follows from the ergodicity
of the action of the group ${\cal M}$ of multiplicators.
\end{proof}

\smallskip\noindent
{\bf Remarks. 1.} This representation may be extended to the group
   $\operatorname{SL}(2,F)$, for this we need to define only one operator, namely, the
   image of the matrix   
$ \left(
\begin{array}{cc}
0&1\\
1&0
\end{array}
\right)$.

{\bf 2.} This representation  of $\operatorname{SL}(2,F)$ 
was firstly  introduced in~\cite{GGV73}
   in a completely different way, but later~\cite{GGV83, GGV85} was interpreted with
   the space $L^2(\hat P_\Gamma)$.

{\bf 3.} This representation is a continual analogue of the classical
   representations of the group $\operatorname{SL}(2,\Bbb R)$ in $L^2(\Bbb R_+)$.
\smallskip

Note that passing from $P_\Ga$ to $\ti P_\Ga$
we do not change the conditional measures $P_\Ga^s$
of $P_\Ga$, given the conic part equal to $s$, 
and modify only
the factor measure on the half-line, 
i.e.~$T'\ti P_\Ga=m_\th\times \operatorname{PD}(\th)\times\bar\nu^\infty$ and
$T\ti P_\Ga=\widetilde{\operatorname{PD}}(\th)\times\bar\nu$, 
where $m_\th$ has density $t^{\th-1}\,/\,\Ga(\th)$, $t>0$
(in particular, $m_1$ is just the Lebesgue
measure on the half-line), 
and $\widetilde{\operatorname{PD}}(\th)=m_\th\times \operatorname{PD}(\th)$.

Recall that the transformation $M_a$ induces a Markovian operator $R_a$ on
the cone $C$ (see Section~\ref{s:PD}).

\begin{corollary}
The $\si$-finite measure $\widetilde{\operatorname{PD}}(\th)$ on the cone $C$
is invariant under the Markovian operator $R_{a}$
for all $a\in\cal M$.
\end{corollary}

\section{Stable processes and general Poisson--Di\-ri\-chlet measures}\label{s:stable}

\begin{definition} \label{def:stable}
Let $\al\in(0,1)$. 
The standard $\al$-stable process
on the space $(X,\nu)$ is
a L\'evy process with L\'evy measure 
\begin{equation}
d\La_\al=\frac{c\al}{\Ga(1-\al)}s^{-\al-1}ds,\qquad s>0,
\label{defstable}
\end{equation}
where $c>0$ is an arbitrary fixed positive number.
\end{definition}

The corresponding infinitely divisible law, i.e.~the 
distribution of the sum of charges with respect to $P_\al$,
is the $\al$-stable law $F_\al$ on $\Bbb R_+$. 

Denote by $P_\al$ the law of the $\al$-stable process.
The Laplace transform of $P_\al$ equals
\begin{equation}
\Bbb E_\al\left[\exp\left(-\int_Xa(x)d\eta(x)\right)\right]=
\exp\left(-c\int_Xa(x)^\al d\nu(x)\right),
\label{stablelapl}
\end{equation}
for an arbitrary measurable function $a:X\to\Bbb R_+$ with
$\int_Xa(x)^\al d\nu(x)<\infty$.

\begin{proposition}[\cite{PY97}]\label{p:pdalpha} 
The simplicial part of an $\al$-stable process with $\al\in(0,1)$ is
the Poisson--Dirichlet distribution $\operatorname{PD}(\al,0)$.
\end{proposition}

\begin{proposition}\label{p:stable}
The conic part of the law of the $\al$-stable process is concentrated on the set
\begin{equation}
\left\{z\in C:\;\lim_{n\to\infty}z_nn^{1/\al}=
\left(\frac{c}{\Ga(1-\al)}\right)^{1/\al}\right\}.
\label{ALasymp}
\end{equation}
\end{proposition}

\begin{proof} In case of an $\al$-stable process, we have
$m(t)=\frac c{\Ga(1-\al)}t^{-\al}$, and~(\ref{ALasymp})
follows immediately from Proposition~\ref{p:asymp}.
\end{proof}

For $s>0$, denote by $\kappa_\al^s$ the 
conditional distribution of the conic part $\kappa_\al$ of the stable process 
on the simplex $\Si_s=\{z=(z_1,z_2,\ldots):\,\sum z_i=s\}$ of monotone sequences
with sum $s$ (i.e.~the conic part of the conditional distribution of the law
$P_\al$ on the set 
$D_s=\{\eta=\sum z_i\de_{x_i}\in D^+:\,\sum z_i=s\}$ 
of positive discrete measures with total charge $s$).

\begin{corollary} 
\label{c:cond}
The homothetic projection of the conditional measure $\kappa_\al^s$ 
on the unit simplex $\Si$
is concentrated on the set
\begin{equation}
\left\{y\in\Si:\;
\lim_{n\to\infty} y_nn^{1/\al}=\frac 1s\left(\frac{c}{\Ga(1-\al)}\right)^{1/\al}\right\}
\label{ALasymp2}
\end{equation}
\end{corollary}

\begin{corollary} \label{c:sing}
The homothetic projections of the measures $\kappa_\al^s$ 
and $\kappa_\al^t$ on the unit simplex $\Si$
are singular for all pairs $s,t>0$, $s\ne t$. Thus the distribution 
$\operatorname{PD}(\al,0)$ is a continual sum of a family of singular distributions.
\end{corollary}

The following statement shows 
how one may recover the conic part of the stable process starting with its simplicial 
part $\operatorname{PD}(\al,0)$. 

\begin{corollary}[\cite{PY97}]\label{c:py}
Let the vector  $Q=(Q_1,Q_2,\ldots)\in\Si$ have the distribution $\operatorname{PD}(\al,0)$. The limit
$L(Q)=\lim_{n\to\infty}n^{1/\al}Q_n$ exists almost surely.
Let $S(Q)=\frac1{L(Q)}\cdot\frac{c}{\Ga(1-\al)}$. Then 
the $\al$-stable process $\eta$ 
on the space $(X,\nu)$ may be represented as
$$
\eta=S(Q)\sum_{i=1}^\infty Q_i\de_{X_i},
$$ 
where $Q$ obeys $\operatorname{PD}(\al,0)$ and $X_1,X_2,\ldots$ is a sequence of 
i.i.d.~variables on $X$
with common distribution $\nu$.
\end{corollary}

The Poisson--Dirichlet distribution $\operatorname{PD}(\al,\th)$ with $\al,\th\ne0$ 
is not the simplicial part of any L\'evy process. However, one may obtain it
as the simplicial part of the process that has density with respect to a stable
process. Namely,
let $\th>-\al$ and consider the law $P_{\al,\th}$ on $D$ which has the density 
\begin{equation}
\frac{dP_{\al,\th}}{dP_\al}(\eta)=\frac{c_{\al,\th}}{\eta(X)^\th}
\label{alth}
\end{equation}
with respect to the $\al$-stable law $P_\al$. Here 
$c_{\al,\th}=c^{\th/\al}\frac{\Ga(\th+1)}{\Ga(\th/\al+1)}$
is a normalizing constant.

\begin{proposition}[\cite{PY97}]
\label{p:althsimp}
The simplicial part of the law $P_{\al,\th}$ is the Poisson--Dirichlet
distribution $\operatorname{PD}(\al,\th)$.
\end{proposition}

\section{The gamma measure as a weak limit of laws of $\al$-stable processes 
when $\alpha$ tends to zero}\label{s:limit}

The purpose of this section is to show that it is natural to consider the gamma process 
as a weak limit of the $\al$-stable processes when $\al\to0$. 
We present several settings of this statement. 

Let $k>0$ and
consider the measure $P_{\al,k}$ on $D$ given by
\begin{equation}
\frac{dP_{\al,k}}{dP_\al}(\eta)=\frac{\exp\left(-\ga\cdot\eta(X)\right)}
{\Bbb E_\al\left[\exp\left(-\ga\cdot\eta(X)\right)\right]}=
e^{\ga^\al}\cdot e^{-\ga\cdot\eta(X)},
\label{dens}
\end{equation}
where $\ga=\frac k{\al^{1/\al}}$.
Denote by $\hat P_{\al,k}$ the law of the process $\ga\cdot\eta$, where 
$\eta$ obeys $P_{\al,k}$. The following theorem was formulated
in~\cite{VY95}. 

\begin{theorem}[\cite{VY95}]\label{th:weak}
The measures $\hat P_{\al,k}$ converge weakly to $P_\Ga$ when $\al\to0$. 
\end{theorem}

\begin{proof}
It follows from~(\ref{dens}),~(\ref{stablelapl}) that
the Laplace transform of the measure $\hat P_{\al,k}$  equals
$$
\Bbb E_{\hat P_{\al,k}}\exp(-f_a(\eta))=
=\exp\left(-\ga^\al\int_X\big((a(x)+1)^\al-1\big)d\nu(x)\right).
$$
But
$$
\ga^\al\big((a(x)+1)^\al-1\big)=
\frac{k^\al}{\al}\big(\al\log(a(x)+1)+o(\al)\big)\to
\log(a(x)+1),
$$ 
as $\al\to0$, hence
$$
\Bbb E_{\hat P_{\al,k}}\exp(-f_a(\eta))\to
\exp\left(-\int_X\log(a(x)+1)d\nu(x)\right)=
\Bbb E_{P_\Ga}\exp(-f_a(\eta)),
$$
and Theorem 7 follows.
\end{proof}

This important result is a key point of the construction 
of two-parameter Poisson--Dirichlet 
distributions $\operatorname{PD}(\al,\th)$ (see Section~\ref{s:gamma}). 
In particular, one obtains the following corollary. 

\begin{corollary}[\cite{PY97}]
For a fixed $\th\ne0$, the distributions $\operatorname{PD}(\al,\th)$ converge to $\operatorname{PD}(0,\th)=\operatorname{PD}(\th)$ 
when $\al\to0$.
\end{corollary}

\begin{corollary}
Let $\ti P_{\al,k}$
 be a measure with constant density $e^{\ga^\al}$
with respect to the $\al$-stable law $P_\al$. 
Denote by $\check P_{\al,k}$ the law of the process
$\ga\cdot\eta$, where $\eta$ obeys $\ti P_{\al,k}$. Then the measures   $\check P_{\al,k}$
converge weakly to the quasi-Lebesgue measure $\ti P_\Ga$.
\end{corollary}

That the quasi-Lebesgue measure is a kind of
the limit case of 
stable processes was suggested in~\cite{Li83}. See
Section~\ref{s:stabletheory} for another understanding of this statement.

There is the following useful way to formalize the transition to a $\si$-finite limit.
Consider the convolution of $n$ copies of the measures $P_{1/n}$ and multiply it by a 
function of $n$. Then the limit in $n$ will be the $\si$-finite measure 
under consideration. In terms of the Laplace transform this procedure is
equivalent to the following relation:
$$
\lim_{n\to\infty}\exp\left(-n\left(x^{1/n}-1\right)\right)=
\frac d{d\al}\exp(-x^\al)\left|_{\al=0}\right.=\frac1x.
$$
Differentiating by $\al$ at the point $\al=0$ is just a ``commutative''
analogue of the main technique used in~\cite{GGV73} for constructing
so-called canonical states and the representation theory of semi-simple
currents for the group $\operatorname{SL}(2,R)$.  

\section{Equivalent definition of stable laws, $\si$--fi\-ni\-te stable measures
and zero--stable laws}
\label{s:stabletheory}

In this section we give another definition of stable laws which is valid for 
$\si$-finite measures.
We describe a piece of theory (to be exposed in details elsewhere) of $\si$-finite 
stable measures showing
that it is natural to consider the Lebesgue measure as
a zero-stable law.

Let $F$ be a distribution on $\Bbb R$.
Consider the distribution $F\times F$ on $\Bbb R\times\Bbb R$. 
Let $\|\cdot\|_\al$ be the $\al$-norm
in the space $(\Bbb R\times\Bbb R)^*$ of linear functionals on $\Bbb R\times\Bbb R$,
i.e.~if $f(x_1,x_2)=a_1x_1+a_2x_2$, then $\|f\|_\al=(|a_1|^\al+|a_2|^\al)^{1/\al}$.
Let us consider only stable laws depending on one parameter $\al\in(0,2]$.
Then the ordinary definition of an $\al$-stable law on $\Bbb R$ is
equivalent to the following one. 

\begin{definition}
The law $F$ is $\al$-stable, if the following condition holds.
If two linear functionals $f_1$ and $f_2$ on $\Bbb R\times\Bbb R$ have the same $\al$-norm,
i.e.~$\|f_1\|_\al=\|f_2\|_\al$, then $f_1$ and $f_2$ have the same distribution
with respect to $F\times F$.
\end{definition}

Note that in case of linear functionals $f_1$ and $f_2$
the equality of distributions is equivalent to 
the existence of a $F\times F$-preserving transformation 
$L:\Bbb R\times\Bbb R\to\Bbb R\times\Bbb R$
such that $f_2=f_1\circ L$. Thus we obtain the following definition of a stable law
which applies to $\si$-finite measures.

\begin{definition}
The measure $F$ (may be $\si$-finite)
is called $\al$-stable, if the following condition holds.
If two linear functionals $f_1$ and $f_2$ on $\Bbb R\times\Bbb R$ have the same $\al$-norm,
i.e.~$\|f_1\|_\al=\|f_2\|_\al$, then there exists a $F\times F$-preserving transformation 
$L:\Bbb R\times\Bbb R\to\Bbb R\times\Bbb R$
such that $f_2=f_1\circ L$.
\end{definition}

Let $a_1a_2\ne0$. If $\al\to0$, then 
\begin{eqnarray*}
2^{-1/\al}\|f\|_\al&=&
\left(\frac{|a_1|^\al+|a_2|^\al}2\right)^{1/\al}=
\left(\frac12|a_1|^\al+\frac12|a_2|^\al\right)^{1/\al}=\\&=&
\left(\frac12+\al\log|a_1|+O(\al^2)+\frac12+\al\log|a_2|+O(\al^2)\right)^{1/\al}=\\&=&
\left(1+\al\log|a_1a_2|+O(\al^2)\right)^{1/\al}\to
\log|a_1a_2|,
\end{eqnarray*}
thus it is natural to consider the quasi-norm $\|f\|_0=|a_1a_2|$ as a natural
limit of $\al$-norms when $\al$ tends to $0$. We obtain the following definition 
of a zero-stable measure on $\Bbb R$.

\begin{definition}
The measure $F$ on $\Bbb R$ is called zero-stable,
if for each two linear functionals $f_1$ and $f_2$ on $\Bbb R\times\Bbb R$
with $\|f_1\|_0=\|f_2\|_0<\infty$, 
there exists a transformation  $L:\Bbb R\times\Bbb R\to\Bbb R\times\Bbb R$
preserving $F\times F$ such that
$f_2=f_1\circ L$. 
\end{definition}

\smallskip\noindent{\bf Remark.}
Note that in case $\alpha=2$, corresponding to the normal law $F$, 
the $F\times F$-preserving transformation $L$, which connects 
functionals of equal 2-norm, is a rotation. In the ``opposite'' case
$\alpha=0$ the corresponding transformation is hyperbolic, of the form
$(a_1,a_2)\mapsto(ca_1,\frac{a_2}c)$.
In both cases we have a linear mapping. Unlike these two extreme cases,
in intermediate cases $\alpha\in(0,2)$ the transformation $L$ is not
linear. 

\begin{proposition}
For any $\beta<1$, the measure on $\Bbb R$ with density
$|x|^{-\beta}dx$ is zero-stable. In particular, 
the Lebesgue measure on $\Bbb R$ is zero-stable. 
\end{proposition}

\begin{proof}
Easy calculation.
\end{proof}

Now we may construct a theory of $\si$-finite stable processes on arbitrary spaces.
Definition~\ref{def:stable} of an $\al$-stable process 
on the space $X$ 
is equivalent to the following one. 
We recall that each bounded Borel function
$a$ on $X$ defines a linear functional $f_a$ on
the space $D$ of finite discrete measures on $X$ by $f_a(\eta)=\int_Xa(x)d\eta(x)$.
Let $\|\cdot\|_\al$ denote the $\al$-norm 
$\|a\|_\al=(\int|a(x)|^\al d\nu(x))^{1/\al}$, and let 
$\|a\|_0=\exp\left(\int_X\log|a(x)|d\nu(x)\right)$. 

\begin{definition}
The law $P_\al$ on $D$ is $\al$-stable, where $\alpha\in[0,2]$,
if for each two linear functionals $f_{a_1}$ and $f_{a_2}$ with 
$\|a_1\|_\al=\|a_2\|_\al<\infty$, 
there exists a transformation $L:D\to D$ preserving the law $P_\al$ such that
$f_{a_2}=f_{a_1}\circ L$. 
\end{definition}

It is easy to check that the quasi-Lebesgue measure $\ti P_\Ga$ satisfies this 
condition. Indeed, if $\|a_1\|_0=\|a_2\|_0$, then $\int_X\log(a_2/a_1)(x)d\nu(x)=0$. 
Hence the multiplicator $M_{a_2/a_1}$ preserves the measure $\ti P_\Ga$ by
Corollary~\ref{c:inv}, and it is obvious that $f_{a_2}=f_{a_1}\circ M_{a_2/a_1}$.
Thus we obtain the following proposition.

\begin{proposition}\label{p:0stable}
The quasi-Lebesgue measure $\ti P_\Ga$ is zero-stable.
\end{proposition}

Note that formulae~(\ref{galapl}) and~(\ref{sigmadens}) imply that
for all $a\in{\cal M}$
\begin{eqnarray*}
&&\Bbb E_{\ti P_\Ga}\left[\exp\left(-\int_Xa(x)d\eta(x)\right)\right]=
\exp\left(-\int_X\log a(x)d\nu(x)\right)\\
&&\qquad=\exp\left(\int_X\log\phi(a(x))d\nu(x)\right),
\end{eqnarray*} 
where $\phi(t)=1/t$ is the Laplace transform of the Lebesgue measure on $\Bbb R_+$
which is zero-stable.
Comparing with~(\ref{laplace}), we see that 
we could define a zero-stable process as a L\'evy process with 
L\'evy measure corresponding to a zero-stable law.

\section{Distributions of linear functionals of the gamma processes and the
Markov--Krein transform}\label{s:mk}

In this section we show that the Markov--Krein identity  
known in the context of Dirichlet processes
may be interpreted
as a formula relating the distribution of a linear functional with
respect to the gamma process and the distribution of the same functional
with respect to the normalized gamma process. This interpretation allows to
prove it immediately, using only a formula for the Laplace
transform of the gamma process.

Given a function $a\in\cal M$, denote by $\mu_a$ the
distribution  
of the linear functional $\eta\mapsto f_a(\eta)=\int_X a(x)d\eta(x)$ 
on $D$ with respect to the law $\bar P_\Ga$ of the normalized gamma process, 
and let $\nu_a$ be the distribution of the function $a$
with respect to the (normalized) parameter measure $\nu$. 
The following property is 
characteristic for the gamma process. 

\begin{theorem}\label{th:mk} 
The measures $\mu_a$ and $\nu_a$ are related by the following integral identity,
\begin{equation}
\int_{\Bbb R}\frac1{(1+zu)^\th}d\mu_a(u)=\exp\left(-\int_X\log(1+zu)^\th d\nu_a(u)\right).
\label{mk}
\end{equation}
\end{theorem}

This formula was first obtained in~\cite{CR90} in the context of Dirichlet processes
by hard analytic arguments (see also simpler proofs in~\cite{DK96} 
and~\cite{KTs98}). But the relation with the gamma process, which is a key
point of our simple proof, has been overlooked. 
Note that the left-hand side of identity~(\ref{mk}) is the generalized 
Cauchy--Stieltjes transform of the distribution $\mu_a$. It is natural to call
the right-hand side the multiplicative version of the generalized 
Cauchy--Stieltjes transform of the distribution $\nu_a$. In view of~(\ref{galapl}), it is equal
to the Laplace transform of the gamma measure $P_\Ga$ calculated on the function $a$.
Thus one may regard formula~(\ref{mk}) as relating an integral transform 
(namely, the Cauchy--Stieltjes transform) of
the distribution $\mu_a$ of the functional $f_a$ with respect to the normalized
gamma process and an integral transform (namely, the Laplace transform) of
its distribution with respect to the non-normalized gamma process. 

In case of $\th=1$, the identity~(\ref{mk}) means that 
the distribution $\mu_a$ is the {\it Markov--Krein transform} of
the measure $\nu_a$. This transform arises in many contexts, such as
the Markov moment problem, continued fractions theory, exponential representations of 
functions of negative imaginary type, the Plancherel growth of Young diagrams, etc.
(see~\cite{Ke98} for a detailed survey).

\begin{proof} Using~(\ref{galapl}), Lemma~\ref{l:ind} 
and the Fubini theorem we obtain that 
the right-hand side of~(\ref{mk}) equals
\begin{eqnarray*}
&&\!\!\!\!\!\!\!\!\!\!\exp\left(-\int_X\log (1+za(x))d\nu(x)\right)=
\Bbb E_{P_\Ga}\left[\exp\left(-z\int_X a(x)d\ga(x)\right)\right]\\
&\qquad=&\Bbb E_{P_\Ga}\left[\exp\left(-z\ga(X)\int_X a(x)d\bar\ga(x)\right)\right]\\
&\qquad=&\Bbb E_{P_{\bar\Ga}}\left[\frac1{\Ga(\th)}\int_0^\infty t^{\th-1}
\exp\left(-t-zt\int_Xa(x)d\bar\ga(x)\right)\right]\\
&\qquad=&\Bbb E_{P_{\bar\Ga}}\left[\frac1{(1+z\int_Xa(x)d\bar\ga(x))^\th}\right],
\end{eqnarray*}
and Theorem follows.
\end{proof}         
    
\smallskip\noindent
{\bf Remarks. 1.} According to a personal communication of P.~Diaconis, the idea
of proving formula~(\ref{mk}) using the Laplace transform formula for the gamma
process was used by F.~Huffer in case of {\it discrete} parameter measure $\nu$ (in
this case the gamma process is just a sum of independent gamma variables and
the normalized gamma process is a random point of a finite-dimensional simplex
obeying a Dirichlet distribution). But the fact that this argument works for
continuous parameter measures as well, which simplifies the proof in essential manner,
seems to have been overlooked.

{\bf 2.} It follows from the known results on the Markov--Krein transform
that the distribution $\mu_a$ of a linear functional $f_a$ is absolutely continuous
(see~\cite{CR90} for an explicit formula for the density).  

{\bf 3.} See~\cite{KTs98} for similar results on the common distributions of
several linear functionals of the Dirichlet process. It is easy to extend the proof
of Theorem~\ref{th:mk} to obtain these results.

\section{The two-parameter generalization of the Markov--Krein formula}
\label{s:genmk}

In \cite{Ts97}
an analogue of the Markov--Krein identity is obtained for the distribution of a linear
functional with respect to the generalized Dirichlet process associated with 
the two-parameter Poisson--Dirichlet distributions $\operatorname{PD}(\al,\th)$. Using the
key idea of Section~\ref{s:mk},
we present here a new proof of this identity based on relation of the
two-parameter Poisson--Dirichlet family to the stable processes.

Let $\al\in(0,1)$ and $\th>-\al$. Denote by $\bar P_{\al,\th}$ the 
normalization (in sense of~(\ref{normaliz})) of the law $P_{\al,\th}$.
(Recall that the simplicial part of this law is $\operatorname{PD}(\al,\th)$, 
see Section~\ref{s:stable}.)
Given an arbitrary measurable function $a:X\to\Bbb R_+$ such that
$\int_Xa(x)^\al d\nu(x)<\infty$,
let $\mu_a$ be the
distribution of the functional $f_a$ with respect to $\bar P_{\al,\th}$, 
and let $\nu_a$ be the distribution of the function $a$
with respect to the (normalized) measure $\bar\nu$. 

\begin{theorem} The measures $\mu_a$ and $\nu_a$ are related by the following
integral identity:

\noindent {\rm 1)} if $\th\ne0$, 
\begin{equation}
\left(\int_{\Bbb R}(1+zu)^{-\th}d\mu_a(u)\right)^{-\frac1\th}=
\left(\int_{\Bbb R}(1+zu)^\al d\nu_a(u)\right)^\frac1\al;
\label{1}
\end{equation}

\noindent {\rm 2)} if $\th=0$,
\begin{equation}
\exp\left(
\int_{\Bbb R}\log(1+zu)^\al d\mu_a(u)\right)=\int_{\Bbb R}(1+zu)^\al d\nu_a(u).
\label{2}
\end{equation}
\end{theorem}

\begin{proof} 1) Denote the left-hand side of the desired identity
by $A^{-1/\th}$ and the right-hand side by $B^{1/\al}$.
Using the identity
$$
\frac1{r^\th}=\frac1{\Ga(\th)}\int_0^\infty t^{\th-1}e^{-rt}dt, 
$$
we obtain
\begin{eqnarray*}
A&=&c_{\al,\th}\Bbb E^\al\left[\left(\eta(X)+z\int_Xa(x)d\eta(x)\right)^{-\th}\right]\\
&=&\frac{c_{\al,\th}}{\Ga(\th)}\Bbb E^\al\left[
\int_0^\infty t^{\th-1}\exp\left(-t\left(\eta(X)+z\int_Xa(x)d\eta(x)\right)\right)dt
\right]\\
&=&\frac{c_{\al,\th}}{\Ga(\th)}\int_0^\infty t^{\th-1}\Bbb E^\al\left[
\exp\left(-\int_Xt(1+za(x))d\eta(x)\right)\right]dt.
\end{eqnarray*}
By the Laplace transform formula~(\ref{stablelapl}), 
the expectation equals precisely $e^{-t^\al B}$, thus
$$
A=\frac{c_{\al,\th}}{\Ga(\th)}\int_0^\infty t^{\th-1}e^{-t^\al B}dt,
$$ 
and~(\ref{1}) follows by changing variables.

2) Follows from~(\ref{1}) by letting $\th\to0$. 
\end{proof}

\section*{Appendix}
\addcontentsline{toc}{section}{Appendix}       
                 
This is a continuation of the discussion on properties of
the  gamma process initiated in the Introduction.
We will start with a new statement which is called the {\it subordinating
property}:

\begin{description}           
\item[6. Subordinating property~\cite{MaGeY00}:] There is the identity in law
$$
(\beta_{\gamma_t},\, t \geq 0)  \stackrel{\operatorname{law}}{=}
\left(\frac{1}{\sqrt2} (\gamma^{(1)}_t - \gamma^{(2)}_t),\,t \geq 0\right),
$$
where the left-hand side, $(\beta_u,\, u \geq 0)$, is a Brownian motion 
independent from a gamma process
$(\gamma_t,\, t \geq 0)$, and in the right-hand side $\gamma^{(1)}$ 
and $\gamma^{(2)}$ are two independent gamma processes.
\end{description}  
                         
Combining the quasi-invariance property~3 (see Introduction)  
and the subordinating property~6 of the gamma process, one easily 
deduces the following 
proposition, which expresses the conditional law of 
$(\gamma_t,\, t \geq 0)$
given $(X_t \equiv \beta_{\gamma_t}, t \geq 0)$, where 
$(\beta_u,\, u \geq 0)$ is a Brownian motion independent 
on $(\gamma_t,\, t \geq 0)$.
           
\begin{proposition}[\cite{MaGeY00}] Conditionally on X, the process 
$(\gamma_t,\,t \geq 0)$ is distributed as $\frac12 T^{(1)}_{V_t}$, $t \geq 0$, where
$(T^{(1)}_a,\, a \geq 0)$ is the process of the first hitting 
times of levels
$a \geq 0$ of a Brownian motion with drift 1 independent of $X$,
and $V_t= \int_0^t |dX_s|$ is a gamma process taken at $(2t)$.          
\end{proposition}
 
The above Proposition suggests the existence of some deep
links between the Brownian motion and the gamma process. We now present
two occurrences of the gamma process related to the Brownian 
motion and to the Bessel processes.
           
\begin{theorem}[\cite{Ber00}] Let $(C_t, t \geq 0)$ denote
the Cauchy process, and for each $x \in \Bbb R$ define
$\mu_x=\sup \{s \leq 1: C_s+xs=\max_{u \leq 1}(C_u+xu)\}$. Then 
$$
(\mu_x,\, x \in \Bbb R)\stackrel{\operatorname{law}}{=} 
(\gamma_{\rho(x)} / \gamma_1,\, x\in \Bbb R),
$$
where $(\gamma_u,\, u \geq 0)$ is the standard gamma process, and
$\rho(x)=\frac12 +\frac{1}{\pi} arctg(x) \equiv P(C_1 \leq x)$.
\end{theorem}                        
                          
See also~\cite[Proposition~1, p.~1535]{Ber99} for some occurrences of
$({\gamma_u}/{\gamma_1},\, u \leq 1)$ in relation with embedded
regenerative sets.

To present the second instance we have in mind, we recall that
the law $Q^{\delta}_0$ on $C(\Bbb R_+,\Bbb R_+)$ of the square 
of a Bessel process of (fractional) dimension $\delta$, 
i.e.~of $\Bbb R_+$-valued diffusion with infinitesimal generator
$2x \frac{d^2}{dx^2}+\delta \frac{d}{dx}$, 
is infinitely divisible, more precisely, 
$Q^{\delta}_0 \* Q^{\delta'}_0=Q^{\delta +\delta'}_0$
for every $\de,\de'\ge0$.
                
In the derivation of the L\'evy--Khinchin representation of
$Q^{\delta}_0$ in~\cite{PY82} was  shown the existence of 
a two-parameter process
$(X^\delta_t,\, \delta \geq 0,\, t \geq 0)$ which may be described as
follows:
\begin{itemize} 
\item[i)] it has homogeneous independent increments in $\delta$, with values
           in the space $C(\Bbb R_+,\Bbb R_+)$;
\item[ii)] for each $\delta > 0$, $(X^\delta_t,\,t \geq 0)$ obeys the law
             $Q^{\delta}_0$;
\item[iii)] for each $t>0$, $(X^{\delta}_t,\, \delta \geq 0)$ is distributed as
              $(2t \gamma_{\delta/2},\, \delta \geq 0)$.
\end{itemize}
          
The limit of the  Bessel processes at the critical value of the
parameter is similar to the situation with stable and gamma processes 
when $\alpha$ tends to zero (see Section~\ref{s:limit}).

\end{document}